\documentclass[a4paper,reqno]{amsart}
\usepackage{amssymb,upref,mathrsfs,url}
\let\mathcal\mathscr
\usepackage[mathscr]{eucal}

\newcommand*{\pd}[2]{\mathchoice{\frac{\partial#1}{\partial#2}}
 {\partial#1/\partial#2}{\partial#1/\partial#2}
 {\partial#1/\partial#2}}

\newcommand{\eval}[2][\right]{\relax
 \ifx#1\right\relax \left.\fi#2#1\rvert}

\newcommand{\envert}[2][\right]{\relax
 \ifx#1\right\relax \left\lvert\else#1\lvert\fi#2#1\rvert}
\let\abs=\envert

\newcommand{\enVert}[2][\right]{\relax
 \ifx#1\right\relax \left\lVert\else#1\lVert\fi#2#1\rVert}
\let\matr=\enVert

\let\kappa\varkappa
\let\phi\varphi

\DeclareMathOperator{\Hom}{Hom}

\DeclareMathOperator{\smbl}{smbl}
\DeclareMathOperator{\ord}{ord}

\DeclareMathOperator{\CDiff}{\mathcal{C}Dif{}f}

\newcommand*{\CDiffskewad}[1]{\CDiff_{(#1)}^{\,\text{\textup{sk-ad}}}}

\newcommand{\cprime}{\/{\mathsurround=0pt$'$}}

\providecommand{\href}[2]{#2}
\newcommand{\urlprefix}{URL }
\newcommand*{\eprint}[2][]{%
\href{http://arXiv.org/abs/#2}{\begingroup \Url{arXiv:#2}}%
}

\DeclareFontFamily{U}{wncyr}{}
\DeclareFontShape{U}{wncyr}{m}{it}{%
  <5><6><7><8><9>gen*wncyi%
  <10><10.95><12><14.4><17.28><20.74><24.88>wncyi10}{}
\DeclareSymbolFont{MathRussLetters}{U}{wncyr}{m}{it}
\DeclareMathSymbol{\re}{\mathalpha}{MathRussLetters}{3}

\newtheorem{theorem}{Theorem}
\newtheorem{proposition}{Proposition}

\theoremstyle{remark}
\newtheorem{remark}{Remark}

\begin{document}

\hfill DIPS-7/2003

\hfill math.DG/0310451

\vspace{1cm}

\title[Integrability conditions]{On the integrability conditions
  for some structures related to evolution differential equations}

\author{P. Kersten}

\address{Paul Kersten \\
  University of Twente \\
  Faculty of Electrical Engineering, Mathematics and Computer Science  \\
  P.O.~Box 217 \\
  7500 AE Enschede \\
  The Netherlands }

\email{kersten@math.utwente.nl}

\author{I. Krasil{\cprime}shchik}

\address{Iosif Krasil{\cprime}shchik \\
  The Diffiety Institute and Independent University of Moscow \\
  B.~Vlasevsky~11 \\
  119002 Moscow \\
  Russia}

\email{josephk@diffiety.ac.ru}
  
\author{A.~Verbovetsky}

\address{Alexander Verbovetsky \\
  Independent University of Moscow \\
  B. Vlasevsky~11 \\
  119002 Moscow \\
  Russia}

\email{verbovet@mccme.ru}
  
\keywords{Evolution equations, Hamiltonian structures, symplectic
  structures, integrability}

\begin{abstract}
  Using the result by D.~Gessler~\cite{Gessler:VSSDSDEq}, we show that
  any invariant variational bivector (resp., variational $2$-form) on
  an evolution equation with nondegenerate right-hand side is
  Hamiltonian (resp., symplectic).
\end{abstract}
\maketitle

\section*{Introduction}

In~\cite{KerstenKrasilshchikVerbovetsky:HOpC}, we described a method
to construct Hamiltonian and symplectic structures on nonlinear
evolution equations. The method was essentially based on the notions
of \emph{variational
  multivector}~\cite{IgoninVerbovetskyVitolo:FLVDOp} and
\emph{variational differential form}. From technical viewpoint, for a
given evolution equation~$\mathcal{E}$, it consisted of two steps:
(1)~solving the linearized equation $\ell_{\mathcal{E}}\phi=0$ in the
so-called \emph{$\ell_{\mathcal{E}}^*$-covering} (resp., the equation
$\ell_{\mathcal{E}}^*\psi=0$ in the
\emph{$\ell_{\mathcal{E}}$-covering}) and (2)~checking the
Hamiltonianity condition $[\![\phi,\phi]\!]=0$, where
$[\![\cdot\,,\cdot]\!]$ denotes the \emph{variational Schouten
  bracket}~\cite{IgoninVerbovetskyVitolo:FLVDOp} (resp., the condition
for~$\psi$ to be symplectic, i.e., closed with respect to a certain
differential in the Vinogradov $\mathcal{C}$-spectral
sequence~\cite{Vinogradov:SSLFCLLTNT}).

Surprisingly enough, it was found out in particular computations that the
second condition always holds true ``by default'' and we still do not
know counterexamples (except for the case of first order equations).
On the other hand, a rather old result by
D.~Gessler~\cite{Gessler:VSSDSDEq} states that all terms
$E_1^{p,n-1}(\mathcal{E})$ of the Vinogradov $\mathcal{C}$-spectral
sequence vanish in the nondegenerate case for $p\ge 3$ (here $n$ is
the number of independent variables).  This fact means exactly that
all variational $2$-forms on nondegenerate evolution equations are
closed and thus symplectic. Since Gessler's proof almost literary
works in the case of multivectors, we immediately obtain that all
bivectors on such equations are Hamiltonian. These facts explain our
experimental results.

We expose the details below. In Section~\ref{sec:gener-jet-bundl},
necessary introduction to the geometry of jet bundles and evolution
equations is presented. Section~\ref{sec:vvf-eq} deals with the
calculus of variational multivectors and forms on evolution equations.
To make our exposition self-contained, we repeat Gessler's proof
from~\cite{Gessler:VSSDSDEq}.  Finally, in Section~\ref{sec:hss} we
derive our main results on the integrability of Hamiltonian and
symplectic structures on nondegenerate evolution equations of
order~$>1$.

\section{Generalities: Jet bundles and evolution equations}
\label{sec:gener-jet-bundl}

Let us fix notation and recall briefly some definitions and results we
will use. For explanations we refer
to~\cite{KrasilshchikVinogradov:SCLDEqMP,KrasilshchikVerbovetsky:HMEqMP,%
  IgoninVerbovetskyVitolo:FLVDOp}.

Let $\pi\colon E\to M$ be a vector bundle over an $n$-dimensional base
manifold~$M$ and $\pi_\infty\colon J^\infty(\pi)\to M$ be the infinite
jet bundle of local sections of the bundle~$\pi$.

In coordinate language, if $x_1,\dots,x_n$,~$u^1,\dots,u^m$ are
coordinates on~$E$ such that $x_i$ are base coordinates and $u^j$ are
fiber ones, then $\pi_\infty\colon J^\infty(\pi)\to M$ is an
infinite-dimensional vector bundle with fiber coordinates
$u^j_\sigma$, where $\sigma=i_1\dots i_{\abs{\sigma}}$ is a symmetric
multi-index.

The basic geometric structure on~$J^\infty(\pi)$ is the \emph{Cartan
  distribution}. In coordinate language, the Cartan distribution is
spanned by the \emph{total derivatives}
\begin{equation*}
  D_i=\pd{}{x_i}+\sum_{j,\sigma}u^j_{\sigma i}\pd{}{u^j_\sigma}.
\end{equation*}

A differential operator on $J^\infty(\pi)$ is called
\emph{$\mathcal{C}$-differential operator} (or \emph{horizontal
  operator}) if it can be written as a sum of compositions of
$C^\infty(J^\infty(\pi))$-linear maps and vector fields that belong to
the Cartan distribution. In coordinates, $\mathcal{C}$-differential
operators are total derivatives operators.

Let $P$ and $Q$ be $C^\infty(J^\infty(\pi))$-modules of sections of
some vector bundles over~$J^\infty(\pi)$. All
$\mathcal{C}$-differential operators from~$P$ to~$Q$ form a
$C^\infty(J^\infty(\pi))$-module denoted by $\CDiff(P,Q)$. More
generally, a map $\Delta\colon P_1\times\dots\times P_k\to Q$ is
called a \emph{multidifferential operator} (of degree~$k$) if it is a
$\mathcal{C}$-differential operator in each argument. Choose elements
$p_i\in P_i$, $i=1,\dots,k$, and consider the operators
\begin{equation*}
  \Delta_i=\Delta(p_1,\dots,p_{i-1},\cdot\,,p_{i+1},\dots,p_k)\colon
  P_i\to Q.
\end{equation*}
Let $l_i$ be the order\footnote{Of course, this number depends on the
  choice of $p$'s; so we define the order as the maximum over all
  possible choices.} of~$\Delta_i$. We define the \emph{symbol}
\begin{equation*}
  \smbl(\Delta)\colon
  S^{l_1}(\Lambda^1(M))\otimes P_1\times\dots
  \times S^{l_k}(\Lambda^1(M))\otimes P_k\to Q,
\end{equation*}
where $S^l$ denotes the symmetric power, of~$\Delta$ as follows. For
any $f\in C^\infty(M)$, let us set
\begin{equation*}
  (\delta_f^{(i)}\Delta)(p_1,\dots,p_k)=
  f\Delta(p_1,\dots,p_k)-\Delta(p_1,\dots,p_{i-1},fp_i,p_{i+1},\dots,p_k)
\end{equation*}
and $\delta_{f_1,\dots,f_l}^{(i)}=
\delta_{f_1}^{(i)}\circ\dots\circ\delta_{f_l}^{(i)}$.  If now
$\omega^i=df_1^i\dots df_{l_i}^i$, $i=1,\dots,k$, are symmetric forms
on~$M$, we set
\begin{equation*}
  (\smbl\Delta)(\omega^1\otimes p_1,\omega^2\otimes p_2,
  \dots,\omega^k\otimes p_k)
  =\delta_{f_1^1,\dots,f_{l_1}^1}^{(1)}\circ\dots
  \circ\delta_{f_1^k,\dots,f_{l_k}^k}^{(k)} (\Delta)(p_1,\dots,p_k).
\end{equation*}
Let $\rho\in J^\infty(\pi)$ and $x=\pi_\infty(\rho)\in M$. Then the
value of the symbol at~$\rho$ is the map
\begin{equation*}
  \eval{\smbl(\Delta)}_\rho\colon
  T_x^*M\otimes P_{1,\rho}\times\dots\times
  T_x^*M\otimes P_{k,\rho}\to Q_\rho
\end{equation*}
polynomially dependent on points of $T_x^*M$ ($P_{i,\rho}$ and
$Q_\rho$ denote here the fibers of the corresponding vector bundles at
the point~$\rho$).

The lift of the de~Rham complex on~$M$ to~$J^\infty(\pi)$ is called
\emph{horizontal de~Rham complex} and is denoted by
\begin{equation*}
  0 \xrightarrow{} C^\infty(J^\infty(\pi)) \xrightarrow{\bar d}
  \bar\Lambda^1(\pi) \xrightarrow{\bar d} \cdots \xrightarrow{\bar d}
  \bar\Lambda^n(\pi)\xrightarrow{}0.
\end{equation*}
The cohomology of the horizontal de~Rham complex are called
\emph{horizontal cohomology} and denoted by $\bar H^q(\pi)$.

The adjoint operator to a $\mathcal{C}$-differential operator
$\Delta\colon P\to Q$ we denote by $\Delta^*\colon \hat Q\to\hat P$,
where $\hat P=\Hom_{C^\infty(J^\infty(\pi))}(P,\bar\Lambda^n(\pi))$.

In coordinates,
\begin{equation*}
  \matr[\bigg]{\sum_\sigma a_{ij}^\sigma D_\sigma}^*
  =\matr[\bigg]{\sum_\sigma(-1)^{\abs{\sigma}}D_\sigma\circ
  a_{ji}^\sigma},
\end{equation*}
where $a_{ij}^\sigma\in C^\infty(J^\infty(\pi))$, and
$D_\sigma=D_{i_1}\circ\dots\circ D_{i_{\abs{\sigma}}}$ for
$\sigma=i_1\dots i_{\abs{\sigma}}$.

Denote by $\CDiffskewad{k}(P,Q)$ the module of $k$-linear
skew-symmetric and skew-adjoint in each argument
$\mathcal{C}$-differential operators $P\times\dots\times P\to Q$.

A $\pi_\infty$-vertical vector filed on~$J^\infty(\pi)$ is called
\emph{evolutionary} if it preserves the Cartan distribution. The Lie
algebra of evolutionary fields is denoted by~$\kappa(\pi)$. It is known
that $\kappa(\pi)$ is naturally isomorphic to the set of sections of
the bundle $\pi_\infty^*(\pi)$; thus $\kappa(\pi)$ is endowed with a
structure of $C^\infty(J^\infty(\pi))$-module.

In local coordinates, the evolutionary field that corresponds to a
section~$\phi=(\phi^1,\dots,\phi^m)$ has the form
\begin{equation*}
  \re_\phi=\sum_{j,\sigma}D_\sigma(\phi^j)\pd{}{u^j_\sigma}.
\end{equation*}

We shall call elements of~$\kappa(\pi)$ \emph{variational vectors}.
Elements of the module $\CDiffskewad{k-1}(\hat\kappa,\kappa)$ will be
called \emph{variational $k$-vector}, while elements of
$\CDiffskewad{k-1}(\kappa,\hat\kappa)$ will be called
\emph{variational $k$-forms}.

One knows that standard constructions and formulas of the calculus of
vector fields and forms on manifolds (the de~Rham differential, inner
product, the Lie derivative, the Schouten bracket) are also valid for
their ``variational'' counterparts, with elements of $\bar H^n(\pi)$
being regarded as ``functions''.

In particular, the Lie derivative on variational vectors
$L_{\re_\phi}\colon\kappa\to\kappa$ takes the form
$L_{\re_\phi}=\re_\phi-\ell_\phi$, where the \emph{linearization
  operator}~$\ell_p$ is defined by the equality
$\ell_p(\alpha)=\re_\alpha(p)$, $\alpha\in\kappa$.

The Lie derivative on variational forms
$L_{\re_\phi}\colon\hat\kappa\to\hat\kappa$ is of the form
$L_{\re_\phi}=\re_{\phi}+\ell_\phi^*$.

The Lie derivative on variational $k$-vectors or $k$-forms satisfies
the equality
\begin{equation*}
  L_{\re_{\phi}}(A)(\xi_1,\dots,\xi_{k-1})
  =L_{\re_{\phi}}(A(\xi_1,\dots,\xi_{k-1}))
  -\sum_iA(\xi_1,\dots,\xi_{i-1},
  L_{\re_{\phi}}(\xi_i),\xi_{i+1},\dots,\xi_{k-1}),
\end{equation*}
where $A$ is a multivector or a form, while $\xi_1,\dots,\xi_{k-1}$
are elements of~$\hat\kappa$ in the former case and elements
of~$\kappa$ in the latter one.

Consider a determined evolution equation \newlength{\DotsLength}
\newlength{\PartLength}
\settowidth{\DotsLength}{$u^1_t=f^1(t,x,u^j_\sigma)$}
\settowidth{\PartLength}{$u^1_t={}$}
\begin{align*}
  &u^1_t=f^1(t,x,u^j_\sigma), \\
  &\hbox to \DotsLength {\dotfill} \\
  &\hbox to \PartLength {\hss $u^m_t={}$}f^m(t,x,u^j_\sigma),
\end{align*}
where $x=(x_1,\dots,x_n)$. We shall interpret it in a geometric way as
the space $\mathcal{E}^\infty=J^\infty(\pi)\times\mathbb{R}$ with the
Cartan distribution generated by the Cartan fields on $J^\infty(\pi)$
and the vector field $D_t=\pd{}{t}+\re_f$, where $t$ is the coordinate
along~$\mathbb{R}$.

The linearization of~$\mathcal{E}^\infty$ is of the form
$\ell_\mathcal{E}=D_t-\ell_f$, while the adjoint linearization is
$\ell_\mathcal{E}^*=-D_t-\ell_f^*$.

Note, that from the above we have
\begin{align}
  \label{eq:3}
  \ell_\mathcal{E}=L_{D_t}&\colon\kappa\to\kappa, \nonumber\\[1ex]
  \ell^*_\mathcal{E}=-L_{D_t}&\colon\hat\kappa\to\hat\kappa.
\end{align}

\section{Variational multivectors and forms on evolution equations}
\label{sec:vvf-eq}

Let $A$ be a (possibly dependent on~$t$) variational multivector or
form on~$J^\infty(\pi)$. If $L_{D_t}(A)=0$ then $A$ is called a
variational multivector or form \emph{on the
  equation~$\mathcal{E}^\infty$}.

\begin{remark}
  Variational bivectors on evolution equations were considered
  in~\cite{KerstenKrasilshchikVerbovetsky:HOpC}.
\end{remark}

\begin{remark}
  From~\eqref{eq:3} it follows that the set of variational $1$-forms
  on~$\mathcal{E}^\infty$ coincides with the term
  $E_1^{1,n-1}(\mathcal{E})=\ker\ell^*_\mathcal{E}$ of the Vinogradov
  spectral sequence.  Similarly, the terms $E_1^{p,n-1}(\mathcal{E})$
  consist of variational $p$-forms.
\end{remark}

The set of variational multivectors and forms on~$\mathcal{E}^\infty$
is closed with respect to all operations that are defined on jet
spaces whenever they are applicable: the differential on variational
forms, inner product, the Schouten bracket, Lie derivative.
\begin{proposition}
  Let $\mathcal{E}^\infty$ be an evolution equation $u_t=f$. For an
  operator $A$ to be a variational $k$-vector or $k$-form
  on~$\mathcal{E}^\infty$ it is necessary and sufficient to have
  \begin{equation}\label{eq:5}
    \nabla(A(\xi_1,\dots,\xi_{k-1}))+\sum_i
    A(\xi_1,\dots,\xi_{i-1},\nabla^*(\xi_i),\xi_{i+1},\dots,\xi_{k-1})=0,
  \end{equation}
  where $\nabla=\ell_\mathcal{E}$ if $A$ is a multivector and
  $\nabla=\ell_\mathcal{E}^*$ if $A$ is a form\textup{;} here
  $\xi_1,\dots,\xi_{k-1}$ are elements of~$\hat\kappa$ in the case of
  multivectors and elements of~$\kappa$ in the case of forms.
\end{proposition}

\begin{proof}
  We have
  \begin{equation*}
    L_{D_t}(A)(\xi_1,\dots,\xi_{k-1}) = L_{D_t}(A(\xi_1,\dots,\xi_{k-1}))
    -\sum_i A(\xi_1,\dots,\xi_{i-1},L_{D_t}(\xi_i),\xi_{i+1}\dots,\xi_{k-1})=0.
  \end{equation*}
  Using~\eqref{eq:3}, we get the result.
\end{proof}

\begin{proposition}
  Let $\mathcal{E}^\infty$ be an evolution equation $u_t=f$. If
  operators $A$,~$A_1$, \dots,~$A_{k-1}$ satisfy the equation
  \begin{equation}
    \label{eq:4}
    \nabla(A(\xi_1,\dots,\xi_{k-1}))+\sum_i
    A_i(\xi_1,\dots,\xi_{i-1},\nabla^*(\xi_i),\xi_{i+1},\dots,\xi_{k-1})=0,
  \end{equation}
  where $\nabla=\ell_\mathcal{E}$ or~$\ell_\mathcal{E}^*$\textup{,}
  then $A_1=A_2=\dots=A_{k-1}=A$.
\end{proposition}

\begin{proof}
  Denote the left-hand side of~\eqref{eq:4} by
  $\Omega(\xi_1,\dots,\xi_{k-1})$.  Then we get
  \begin{equation*}
    \Omega(\xi_1,\dots,\xi_{i-1},t\xi_i,\xi_{i+1}\dots,\xi_{k-1})
    -t\Omega(\xi_1,\dots,\xi_{k-1})
    =\pm (A(\xi_1,\dots,\xi_{k-1})-A_i(\xi_1,\dots,\xi_{k-1}))=0.
  \end{equation*}
\end{proof}

\begin{remark}
  The last proposition shows that computing variational multivectors
  and forms on an equation amounts to solving equation $\nabla(s)=0$
  on the $\nabla^*$-covering
  (see~\cite{KerstenKrasilshchikVerbovetsky:HOpC} for the definition
  of $\Delta$-coverings associated to a $\mathcal{C}$-differential
  operator $\Delta$).
\end{remark}

\begin{theorem}
  \label{sec:vari-mult-forms-2}
  Suppose that the symbol of the $\mathcal{C}$-differential
  operator~$\ell_f$ is nonsingular on a dense open subset
  of~$\mathcal{E}^\infty$ and the order of the operator~$\ell_f$ is
  greater than~$1$.  Then there are no nonzero operators~$A$ that
  satisfy equation~\eqref{eq:5} for~$k\geq3$.
\end{theorem}

\begin{proof}[Proof~\textup{({\cite[Th.~$3$]{Gessler:VSSDSDEq}})}]
  Equation~\eqref{eq:5} can be written in the form
  \begin{equation}
    \label{eq:15}
    \pm D_t(A)(\xi_1,\dots,\xi_{k-1})
    +\nabla'(A(\xi_1,\dots,\xi_{k-1}))+\sum_i
    A(\xi_1,\dots,\xi_{i-1},\nabla'^*(\xi_i),\xi_{i+1},\dots,\xi_{k-1})=0,
  \end{equation}
  where $\nabla'=\ell_f$ if $A$ is a multivector and
  $\nabla'=\ell_f^*$ if $A$ is a form. Take a point
  $\rho\in\mathcal{E}^\infty$ such that the symbol
  $\lambda=\eval{\smbl(\nabla')}_\rho$ is nondegenerate at~$\rho$.
  Let $\theta=\sum_{i=1}^n\theta_i\eval{dx_i}_\rho$ be a covector,
  so that, in coordinates,~$\lambda$ is an $m\times m$ matrix
  $\lambda=\matr[\big]{\lambda^i_j}$, where $\lambda^i_j$ are
  homogeneous polynomials in $\theta$'s of degree $l=\ord(\nabla')$.
  Denote the components of the symbol $a=\eval{\smbl(A)}_\rho$ by
  $a^j_{i_1\dots i_{k-1}}(\theta^1,\dots,\theta^{k-1})$,
  $\theta^p=(\theta^p_1,\dots,\theta^p_n)$. Then the symbol of
  equation~\eqref{eq:15} takes the form
  \begin{equation}
    \label{eq:6}
    \sum_{j=1}^m\lambda_j^i(\theta^1+\dots+\theta^{k-1})a^j_{i_1\dots i_{k-1}}
    +(-1)^l\sum_{p=1}^{k-1}\sum_{j=1}^m a^i_{i_1\dots
      i_{p-1}ji_{p+1}i_{k-1}}
    \lambda_j^{i_p}(\theta^p)=0,
  \end{equation}
  where $1\leq i$, $i_1$, \dots,~$i_{k-1}\leq m$.

  System~\eqref{eq:6} can be considered as a linear system of
  algebraic equations with polynomial coefficients over~$\mathbb{C}$.
  Let us show that the determinant of this system does not vanish.

  Since $\lambda=\lambda(\theta)$ is nonsingular, there exists
  $v\in\mathbb{C}^m$ such that $\det\lambda(v)\neq 0$. One can assume
  that $\lambda(v)$ has an upper triangular form, $\lambda_j^i(v)=0$
  if $i\geq j$ and $\lambda^i_i(v)\neq 0$. Then for any
  $\alpha\in\mathbb{C}$ the matrix $\lambda(\alpha
  v)=\alpha^l\lambda(v)$ has also an upper triangular form. Since
  $l=\ord\ell_f\geq 2$ and $k\geq 3$, there exist
  $\alpha_p\in\mathbb{C}$, $p=1,\dots,k-1$, such that for any $1\leq
  i$,~$i_1$, \dots,~$i_{k-1}\leq m$
  \begin{equation}
  A_{ii_1\dots i_{k-1}}=\lambda^i_i(v)
  (\alpha_1+\dots+\alpha_{k-1})^l+(-1)^l
  \sum_{p=1}^{k-1}\lambda_{i_p}^{i_p}(\alpha_p)^l\neq 0.
\end{equation}
Put $\theta^i=\alpha_i v$. Then system~\eqref{eq:15} is upper
triangular with respect to the lexicographic order of indexes, with
diagonal entries $A_{ii_1\dots i_{k-1}}\neq0$. Hence, the determinant
of system~\eqref{eq:15} does not equal to zero, thus $a=0$. Therefore,
the symbol of~$A$ vanishes on a dense subset of~$\mathcal{E}^\infty$,
so that~$A=0$.
\end{proof}

\section{Hamiltonian and symplectic structures}
\label{sec:hss}
Recall that a bivector~$A$ on equation~$\mathcal{E}^\infty$ is said to
be \emph{Hamiltonian} if $[\![A,A]\!]=0$, where
$[\![\cdot\,,\cdot]\!]$ is the \emph{variational Schouten
  bracket}~\cite{IgoninVerbovetskyVitolo:FLVDOp,KerstenKrasilshchikVerbovetsky:HOpC}.
Two structures~$A$ and~$B$ are \emph{compatible} (or constitute a
\emph{pencil}) if $[\![A,B]\!]=0$. Respectively, a \emph{symplectic
  structure} on~$\mathcal{E}^\infty$ is a closed variational form.

Now, our main result is obtained by reformulating the previous
theorem.
\begin{theorem}
  Assume that $\mathcal{E}=\{u_t=f\}$ is an evolution equation and
  $\ell_f$ satisfies the hypothesis of
  Theorem~\ref{sec:vari-mult-forms-2}. Then any variational bivector
  on $\mathcal{E}$ is Hamiltonian\textup{,} any two Hamiltonian
  structures are compatible and any variational $2$-form is
  symplectic.
\end{theorem}
\begin{proof}
  Indeed, let~$A$ and~$B$ be bivectors. Then $[\![A,B]\!]$ is a
  $3$-vector and thus vanishes.  Similarly, the differential of any
  variational $2$-form is a $3$-form and therefore equals zero.
\end{proof}

\begin{remark}
  The hypothesis of Theorem~\ref{sec:vari-mult-forms-2} comprises two
  conditions: (1)~the order of~$\mathcal{E}$ is to be~$>1$; (2)~the
  symbol of~$\ell_f$ is to be nondegenerate. The first one is really
  essential. As for the second condition, it seems that it may be
  weakened. At least in some computation (e.g., for the Boussinesq
  equation, see~\cite{KerstenKrasilshchikVerbovetsky:HOpC}) all
  bivectors are automatically Hamiltonian.
\end{remark}

\begin{remark}
  The proof of Theorem~\ref{sec:vari-mult-forms-2} does not use the
  fact that the operator~$A$ is skew-symmetric and holds also for
  symmetric $\mathcal{C}$-differential operators. This means that
  equations satisfying the hypothesis of the theorem may admit
  \emph{linear} Hamiltonian and symplectic structures only.
\end{remark}

\end{document}